\documentclass{iopart}
\usepackage{iopams}
\pdfoutput=1
\usepackage{latexsym,color,wrapfig}
\usepackage{diagbox}
\usepackage{comment}
\usepackage{amsfonts}
\usepackage{graphicx}
\usepackage{subfig}
\usepackage{url}
\def\sqr#1#2{$\vcenter{\hrule height.#2pt
\hbox{\vrule width.#2pt height#1pt \kern#1pt \vrule width.#2pt} \hrule
height.#2pt}$}
\def\qed{\sqr53}

\def\p{\prime}
\def\G{\tilde{G}}

\newtheorem{lemma}{Lemma}
\newtheorem{theorem}{Theorem}

\newtheorem{definition}{Definition}

\newtheorem{remark}{Remark}
\newtheorem{example}[theorem]{Example}

\begin{document}
\title[Writhe-like invariants of alternating links]{Writhe-like invariants of alternating links}
\author{Yuanan Diao$^{\ast}$$^{\ddag}$, Van Pham$^{\ddag}$}
\address{
$^{\ddag}$ Department of Mathematics and Statistics\\
University of North Carolina at Charlotte\\
 Charlotte, NC 28223, USA\\
$^{\ast}$ Correspondance E-mail: ydiao@uncc.edu}

\begin{abstract}
It is known that the writhe calculated from any reduced alternating link diagram of the same (alternating) link has the same value. That is, it is a link invariant if we restrict ourselves to reduced alternating link diagrams. This is due to the fact that reduced alternating link diagrams of the same link are obtainable from each other via flypes and flypes do not change writhe. In this paper, we introduce several quantities that are derived from Seifert graphs of reduced alternating link diagrams. We prove that they are ``writhe-like" invariants in the sense that they are also link invariants among reduced alternating link diagrams. The determination of these invariants are elementary and non-recursive so they are easy to calculate. We demonstrate that many different alternating links can be easily distinguished by these new invariants, even for large, complicated knots for which other invariants such as the Jones polynomial are hard to compute. As an application, we also derive an if and only if condition for a strongly invertible rational link.
\end{abstract}

\section{Introduction}

\begin{figure}[htb!]
\begin{center}
\includegraphics[scale=0.5]{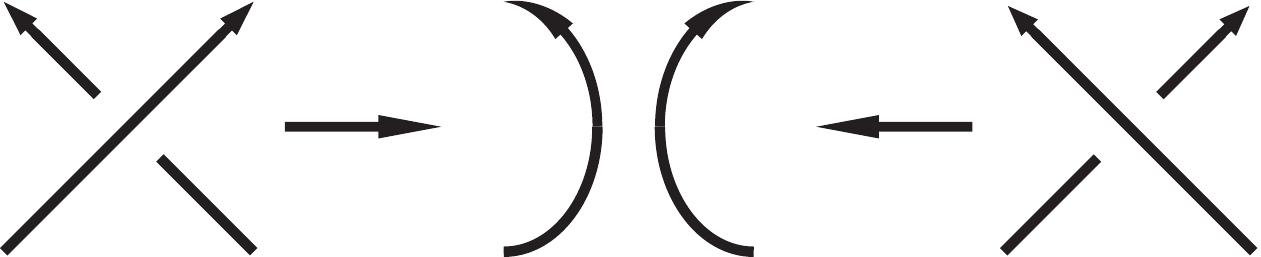}
\end{center}
\caption{How a crossing in an oriented link diagram is smoothed. The crossing on the left is a positive crossing and the one on the right is a negative crossing.
}
\label{smoothing}
\end{figure}

In this paper, the authors wish to derive new link invariants that are ``writhe-like" in the sense that they are link invariants if we restrict ourselves to reduced alternating link diagrams. More specifically, these invariants will be defined using quantities extracted from the Seifert graphs of alternating link diagrams. 
Let $D$ be a link diagram of an oriented (alternating) link $L$. If we {\em smooth} each crossing of $D$ as shown in Figure \ref{smoothing}, then we obtain a set of disjoint topological circles, the collection of which (together with the information on the crossings being smoothed) is called the {\em Seifert circle decomposition} (or $s$-decomposition for short) of $D$. In general, there can be multiple crossings between two Seifert circles of $D$ and these crossings can have different signs. However, it is a well known fact that alternating link diagrams are homogeneous, that is, the crossings on one side of any $s$-circle are of the same sign, while the crossings on the other side of the $s$-circle are of the opposite sign. It follows that all crossings between two Seifert circles of $D$ have the same sign if $D$ is an alternating link diagram.

\begin{definition}{\em
Let $D$ be an oriented link diagram and $S(D)$ be its $s$-decomposition. We construct a graph $G_S(D)$ from $S(D)$ by identifying each Seifert circle ($s$-circle) to a vertex of $G_S(D)$. If there exist $k\ge 1$ crossings between two Seifert circles $C_1$ and $C_2$ in $S(D)$, then the two corresponding vertices  $v_1$ and $v_2$ are connected by $k$ edges. Otherwise, there does not exist an edge between the two vertices. $G_S(D)$ is called the {\em Seifert graph} of $D$.}
\label{sgraph}
\end{definition}

\begin{figure}[htb!]
\begin{center}
\includegraphics[scale=.9]{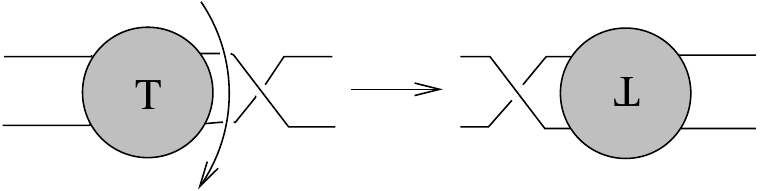}
\end{center}
\caption{A flype on the tangle $T$ moves the crossing $O$ from one side of $T$ to the other, while turning $T$ by $180^\circ$ around the horizontal axis in the direction as shown.
}
\label{flype}
\end{figure}

\medskip
A {\em flype} is a move that can be applied to a link diagram as shown in Figure \ref{flype}. 
The well known Tait's Flype Conjecture states that if $D_1$ and $D_2$ are reduced alternating diagrams, then $D_2$ can be obtained from $D_1$ by a sequence of flypes if and only if $D_1$ and $D_2$ represent the same link. This conjecture was proved to be true by Menasco and Thistlethwaite \cite{MT1,MT2} in 1990. 

\medskip
The writhe of a link diagram is defined as the total number of positive crossings minus the total number of negative crossings of a link diagram. Since a flype does not change the writhe, the result of Menasco and Thistlethwaite \cite{MT1,MT2} implies that writhe is a link invariant if we restrict ourselves to reduced alternating link diagrams. In general the writhe of a link diagram is not an invariant in the equivalence class induced by all Reidemeister moves, because a type I Reidemeister move changes the writhe by $+1$ or $-1$. It is neither an invariant even if we restrict ourselves to minimal link diagrams. A well known example is the Perko pair consisting of knots $10_{161}$ and $10_{162}$ enumerated in Rolfsen's table \cite{Rolfsen}. Of course these knots are non-alternating. 

\medskip
We say that a quantity derived from a link diagram is a {\em writhe-like invariant} if, like the writhe, it is a link invariant when restricted to reduced alternating link diagrams. The main goal of this paper is to introduce several writhe-like invariants. In Section \ref{graphbasic}, we introduce some specific concepts/terminology in graph theory which are used in our paper. In Section \ref{sgraphs}, we examine in detail how flypes change the Seifert graphs of reduced alternating link diagrams. We state our main result in Section \ref{main} where we introduce a set of writhe-like invariants. Finally, in Section \ref{applications}, we show a few applications of some of these writhe-like invariants. 

\section{Background Knowledge in Graph Theory}\label{graphbasic}

We will assume that our reader has the basic knowledge in graph theory and is familiar with the common graph theory concepts and terminology such as vertices, edges, faces, cycles, valences of vertices. The definitions of these can be found in any standard graph theory textbook, for example in \cite{West}.

\medskip
In this section we will introduce certain concepts and terminology that are more specific to this paper. Our graphs are those derived from Seifert circle decompositions of (oriented) reduced alternating link diagrams which will be introduced in the next section. As such, the graphs discussed here are all plane graphs and are bipartite (hence it does not contain any loop edges, nor cycles of odd length). Also, a graph in this paper is not necessarily a simple graph, that is, multiple edges between two vertices are possible. If there is only one edge connecting two vertices, then the edge will be called a {\em single edge}, otherwise edges connecting two vertices are called {\em multiple edges} if there is a need to clarify this.  At times it is more convenient to treat such a graph $G$ as an edge weighted simple graph $\G$, called the {\em simple representation} of $G$, by identifying all edges incident to the same two vertices $v_1$ and $v_2$ so that the weight of the edge in $\G$ connecting $v_1$ and $v_2$ equals the number of edges in $G$ incident to $v_1$ and $v_2$. Notice that if $G$ is connected and bipartite, then so is $\G$. Furthermore, if the edges between two vertices $v_1$ and $v_2$ form an edge cut for $G$, then these edges corresponds to a bridge edge in $\G$. Let $G$ be a graph and $H$ be a subgraph of $G$, that is, the vertex set and edge set of $H$ are subsets of the vertex set and edge set of $G$ respectively. We say that $H$ is an {\em induced subgraph} of $G$ if $H$ has the property that if two vertices $v_1$, $v_2$ are in $H$, then all the edges of $G$ connecting $v_1$ and $v_2$ also belong to the edge set of $H$. If $H$ is an induced subgraph of $G$, we use the notation $G\setminus H$ to denote the subgraph of $G$ whose vertex set and edge set are $(V(G)\setminus V(H))\cup (V(G)\cap V(H)))$ and $E(G)\setminus E(H)$ respectively.

\medskip
\begin{definition}\label{block_definition}{\em
Let $B$ be an induced and connected subgraph of $G$. We say that $B$ is a {\em block} of $G$  
if $B$ contains no cut vertices and is maximal in the sense that if there is another induced and connected subgraph $B^\p$ of $G$ that contains no cut vertices, and $B$ is a subgraph of $B^\p$, then $B=B^\p$.
}
\end{definition}

\medskip
\begin{lemma}\label{block_types}
Let $B$ be a block of a connected graph $G$ with at least one edge, then $\tilde{B}$ is either a single edge graph and this edge is a bridge edge of $\G$ or every edge of $\tilde{B}$ belongs to a cycle of even length at least 4 in $\tilde{B}$.
\end{lemma}

\medskip
\begin{proof}
If $B$ is a single vertex then it is not a block since it belongs to any block that contains an edge connected to it hence it does not satisfy the maximal condition. If $\tilde{B}$ is a single edge graph with vertices $v_1$ and $v_2$ and this edge is a bridge edge of $\G$ then $B$ must be a block. To see this we only need to show that it is maximal since it obviously does not contain any cut vertices of its own. If it is not maximal, then there exists a block $B_1$ that contains it as a proper subgraph. This is impossible since at least one of $v_1$ and $v_2$ would be a cut vertex of $B_1$. That is, (i) and (ii) are indeed possible candidates for $B$. If $\tilde{B}$ is not a single edge graph then let us consider an edge $e\in E(\tilde{B})$. Let $v_1$, $v_2$ be the end vertices of $e$. Notice that $\tilde{B}$ is a simple graph (hence it does not contain any cycles of length 2) and it is still free of cut vertices. By the assumption,  $\tilde{B}$ has at least one more vertex $v_3$. It is a well known fact that a simple graph with 3 or more vertices is free of cut vertices if and only if every edge in it belongs to a cycle. Thus $e$ belongs a  cycle in $\tilde{B}$ with an even length at least 4 since $\tilde{B}$ is also bipartite. \qed
\end{proof}

\medskip
\begin{definition}\label{kcut}{\em
A pair of vertices $\{v_1,v_2\}\subset V(G)$ is called a {\em 2-cut} if $v_1$ and $v_2$ belong to a block $H$ of $G$, and deleting $v_1$ and $v_2$ (and edges incident to them) from $H$ results in a disconnected graph.  
}
\end{definition}

\begin{definition}{\em
Let $G$ with an induced and connected subgraph $H$. If $v$ is a cut vertex, $v\in V(H)$ and $H\cap (G\setminus H)=\{v\}$, then $H$ is said to be $\{v\}$-dependent. Similarly, if $\{v_1,v_2\}$ is a 2-cut and $\{v_1,v_2\}\subset V(H)$ such that $H\cap (G\setminus H)=\{v_1,v_2\}$, then $H$ is said to be $\{v_1,v_2\}$-dependent.
}
\end{definition}

\medskip
\begin{remark}{\em
Let $H$ be a block in $G$ and $\{v_1,v_2\}\subset V(H)$ is a 2-cut of $H$ such that deleting $v_1$ and $v_2$ (and edges incident to them) from $H$ results in a disconnected graph, then $H$ must contain at least one more vertex other than $v_1$ and $v_2$, since otherwise deleting $v_1$, $v_2$ and the edges incident to them from $H$ will result in an empty graph which is not a disconnected graph. By Lemma \ref {block_types}, $H$ must contain at least 4 vertices including $v_1$ and $v_2$ since $\tilde{H}$ must contain cycles of length at least 4.
}
\end{remark}

\medskip
\begin{definition}\label{necklace_definition}{\em
Let $H$ be a block of $G$. We say that $H$ admits a {\em necklace decomposition} if there exist a set of vertices $\{v_1, v_2, ..., v_k\}\subset V(H)$ ($k\ge 2$) and induced subgraphs $H_1$, ..., $H_{k}$ of $H$ such that (i) $E(H_i)\cap E(H_j)=\emptyset$ if $i\not=j$, $1\le i, j\le k$; (ii) $E(H)=\cup_{1\le j\le k}E(H_j)$; (iii) $V(H_i)\cap V(H_{i+1})=\{v_{i+1}\}$ for $1\le i\le k-1$, $V(H_k)\cap V(H_{1})=\{v_{1}\}$, and $V(H_i)\cap V(H_{j})=\emptyset$ for all other $H_i$, $H_j$ pairs; (iv) The set $\{v_1, v_2, ..., v_k\}$ is maximal in the sense that if there exist another set of vertices $\{v^\p_1, v^\p_2, ..., v^\p_{k^\p}\}\subset V(H)$ and induced subgraphs $H^\p_1$, ..., $H^\p_{k^\p}$ of $H$ satisfying conditions (i) to (iii), and  $\{v_1, v_2, ..., v_k\}\subset\{v^\p_1, v^\p_2, ..., v^\p_{k^\p}\}$, then we must have $\{v_1, v_2, ..., v_k\}=\{v^\p_1, v^\p_2, ..., v^\p_{k^\p}\}$ (and consequently  $k^\p=k$, $H_j=H^\p_j$ for all $j$).
}
\end{definition}

\begin{figure}[htb!]
\begin{center}
\includegraphics[scale=0.6]{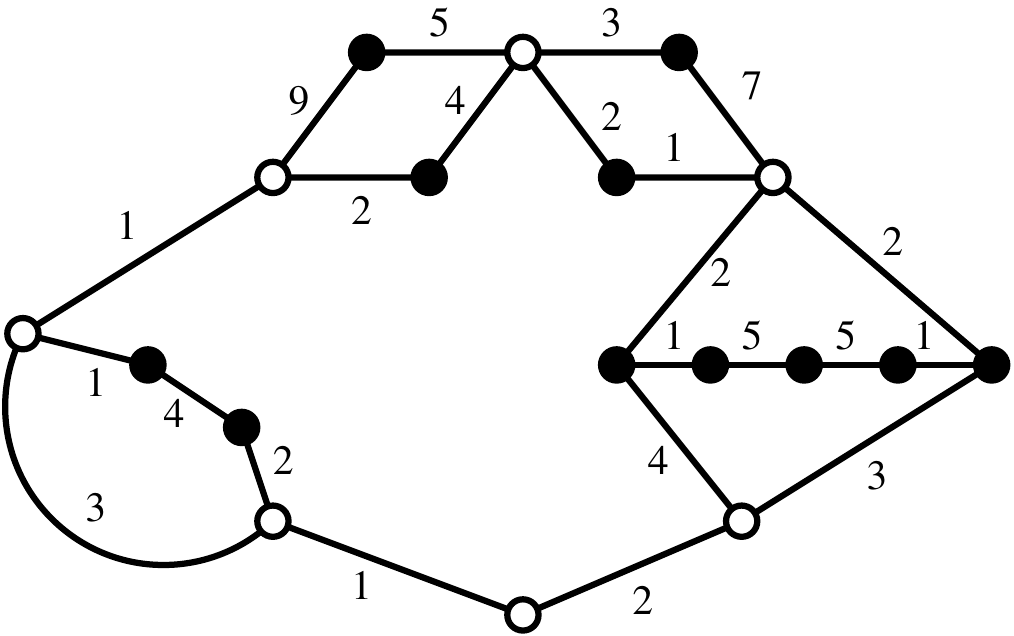}
\end{center}
\caption{A singular necklace that contains 7 beads (whose terminal vertices are marked by circles) with 2 of which being single edge graphs. Multiple edges between two vertices are represented by a single edge with a weight number indicating the number of these edges.
}
\label{chainedbeads}
\end{figure}

We shall call each $H_j$ in Definition \ref{necklace_definition} a {\em bead} (of the necklace decomposition), and the two vertices $v_j$, $v_{j+1}$ ($v_k$, $v_1$ in the case that $j=k$) the {\em terminals} of the bead $H_j$. Notice that a bead can be a single edge graph. We will say that a necklace decomposition is {\em singular} if the decomposition contains at least three beads and at least one bead is a single edge graph. Figure \ref{chainedbeads} is an example of a necklace decomposition of a graph block. We would like to point out that a block may admits more than one necklace decomposition.

In \cite{Whitney}, Whitney introduced a graph operation called {\em Whitney flip}. In this paper we will use several similar operations, which are defined below.

\begin{figure}[htb!]
\begin{center}
\includegraphics[scale=0.7]{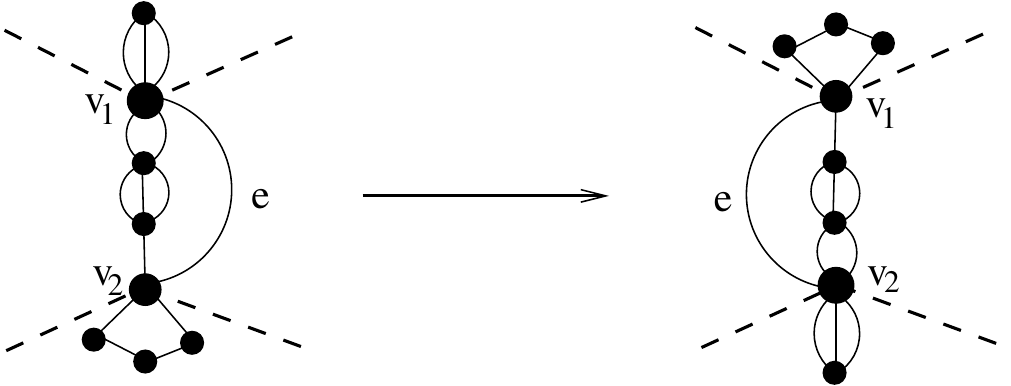}
\end{center}
\caption{An example of Type A Whitney flip.}
\label{WhitneyFlipA}
\end{figure}

\medskip
\begin{definition}\label{TypeAFlip}{\em
Let $\{v_1,v_2\}$ be a 2-cut of a plane graph $G$ and let $H$ be an induced and connected subgraph of $G$ that is $\{v_1,v_2\}$-dependent.  Let $H^\p$ be the mirror image of $H$ with $v_1^\p$ and $v_2^\p$ being its vertices corresponding to $v_1$ and $v_2$. In the case that $v_1$ and $v_2$ are connected by an edge $e$, then a {\em Type A Whitney flip}  is the operation that attaches $H^\p$ to $G\setminus H$ by identifying $v_1^\p$ to $v_2$ and $v_2^\p$ to $v_1$ as shown in Figure \ref{WhitneyFlipA}. 
}
\end{definition}

\begin{figure}[htb!]
\begin{center}
\includegraphics[scale=0.7]{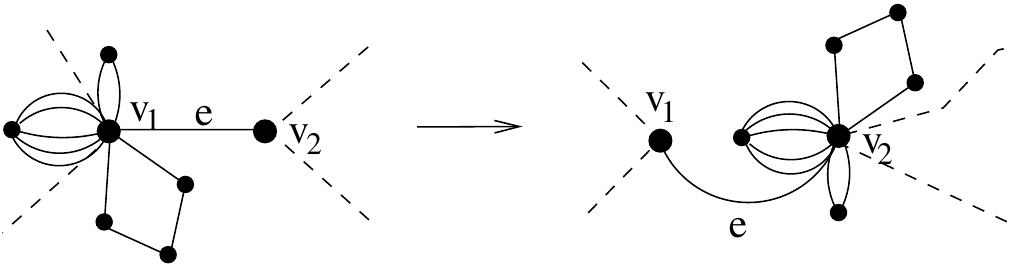}
\end{center}
\caption{An example of Type B Whitney Flip.
}
\label{WhitneyFlipB}
\end{figure}

\medskip
\begin{definition}\label{TypeBFlip}{\em
Let $v_1$, $v_2$ be two vertices of a plane graph $G$ that are connected by a single edge $e$. Let $H$ be an induced and connected subgraph of $G$ that is $v_1$-dependent. The operation that attaches its mirror image $H^\p$ to $G\setminus H$ by identifying $v_1^\p$ (the vertex in $H^\p$ corresponding to $v_1$) to
$v_2$ while keeping $e$ intact is called a {\em Type B Whitney flip}. An example is shown in Figure \ref{WhitneyFlipB}.
}
\end{definition}

\begin{figure}[htb!]
\begin{center}
\includegraphics[scale=0.7]{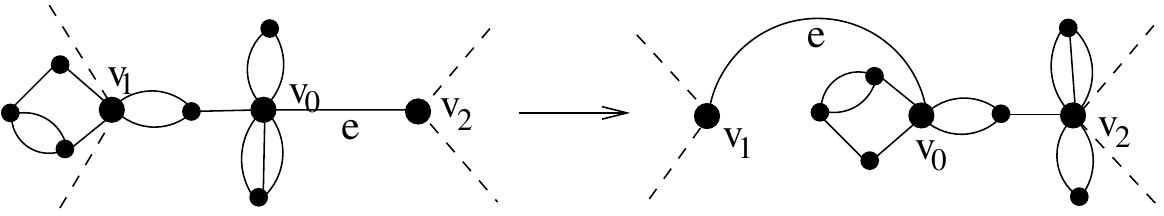}
\end{center}
\caption{An example of Type C Whitney Flip.
}
\label{WhitneyFlipC}
\end{figure}

\medskip
\begin{definition}\label{TypeCFlip}{\em
Let $v_0$, $v_1$, $v_2$ be three vertices of a plane graph $G$ such that there is a single edge $e$ between $v_0$, $v_2$ and $v_1$ is a 1-cut of $G-e$. Let $H$ be an induced and connected subgraph of $G$ that is $v_1$-dependent in $G-e$ and let $H^\p$ be the mirror image of $H$ with $v_0^\p$, $v_1^\p$ being its vertices corresponding to $v_0$ and $v_1$. The operation that attaches $H^\p$ to $G\setminus H$ by identifying $v_0^\p$ to $v_2$, $v_1^\p$ to $v_0$, and moves the edge $e$ to be between $v_0$ and $v_1$ is called a {\em Type C Whitney flip}. An example is shown in Figure \ref{WhitneyFlipC}.
}
\end{definition}

\medskip
\begin{remark}\label{Whitney_flip_remark}{\em 
For the purpose of this paper, we will only consider graphs without single bridge edges. Thus we are only interested in what the Whitney flips defined above can do to such a graph. We observe that a  Whitney flip of Type A or B can only move a block of $G$ attached to one cut vertex to another cut vertex without changing the internal structure of the block. On the other hand, a Whitney flip of Type C can only be applied to a block $B$ that has a singular necklace decomposition and must use a bead $H^\p$ that is a single edge graph as the edge $e$ in the definition of the flip. The flip exchanges the position of $e$ and a chain of beads (one or more consecutive beads, that is). Some blocks attached to cut vertices that belong to either this chain of beads or $H^\p$ may be moved to another terminal without internal structure change as the by-product of this flip without any internal structural change.
}
\end{remark}

\medskip
\begin{definition}\label{flipclass}{\em
Let $G_1$ and $G_2$ be two plane and bipartite graphs. If $G_1$ and $G_2$ are related by a number of Whitney flips of any types A to D, then we say that they belong to the same {\em flip equivalence class} and we write this as $G_1 \stackrel{f}{\equiv} G_2$.}
\end{definition}

\medskip
\begin{remark}\label{Whitney_flip_remark2}{\em 
Let $G_1$ and $G_2$ be two plane and bipartite graphs such that $G_1 \stackrel{f}{\equiv} G_2$. By Remark \ref{Whitney_flip_remark} and Lemma \ref{block_types}, if $\{e_1,e_2,...,e_m\}$ is the set of edges between two vertices in $G_1$, then they will remain the set of edges between two vertices in $G_2$. That is, an edge of weight $m$ in $\G_1$ will remain an edge of weight $m$ in $\G_2$. Furthermore, a cycle in $\G_1$ will remain a cycle in $\G_2$ with the same set of edges (only the cyclic order of the edges along the cycle may differ).
}
\end{remark}

\medskip

\section{Seifert Graphs of Reduced Alternating Link Diagrams }\label{sgraphs}

\medskip
\begin{definition}\label{level}{\em
Let $D$ be a connected diagram of an oriented link and consider $S(D)$. If an $s$-circle $C$ is not bounded inside of any other $s$-circle, then we say that $C$ has {\em level number} 0 and we write this as $\ell(C)=0$. Apparently for any diagram $D$, $S(D)$ has at least one level 0 $s$-circle. In other words, $G_S(D)$ has at least one root. If an $s$-circle $C$ is bounded within a level 0 $s$-circle, but not bounded inside any other $s$-circles, then we say that it has level number $-1$ and we write this as $\ell(C)=-1$. Repeating this process, the level number of any $s$-circle can be defined recursively: assume that all $s$-circles of level $\ge -k$ have been defined, then an $s$-circle $C$ that is bounded within a level $-k$ $s$-circle (but not other $s$-circles that have not been defined at this point) is assigned a level number $-k-1$: $\ell(C)=-k-1$. 
}
\end{definition}

\medskip
\begin{definition}\label{s_block}{\em
Consider the case when a level $-k+1$ $s$-circle bounds one or more level $-k$ $s$-circles. Since the original diagram is connected, the diagram corresponding to these $s$-circles (together with the crossings they share among them) is also connected. We call the collection of these $s$-circles an {$s$-block} and define $-k$ as the  {\em level number} of the $s$-block. The collection of  level 0 $s$-circles form the level $0$ $s$-block. In particular, if there is only one level 0 $s$-circle, then the level $0$-block contains only one $s$-circle. By an abuse of terms, we will use the same name for the counterpart of an $s$-block in $G_S(D)$.
}
\end{definition}

\medskip
\begin{remark}\label{sign_remark}{\em 
It is not hard to prove that $G_S(D)$ can be realized as a plane graph in general (so we will assume that $G_S(D)$ is a plane graph in the rest of the paper), and $G_S(D)$ is bipartite. Furthermore, in the case that $D$ is a reduced alternating link diagram, it is well known that $D$ is homogeneous, that is, the crossings on one side of any $s$-circle are of the same sign, while the crossings on the other side of the $s$-circle are of the opposite sign. It follows that all crossings belonging to the same $s$-block have the same sign, and $s$-blocks whose level numbers differ by an odd integer contain crossings of different signs.
}
\end{remark}

\medskip
\begin{remark}\label{s_block_remark}{\em 
If an $s$-block is bounded by a Seifert circle $C$, then the Seifert graph $H$ of this $s$-block is an induced subgraph of $G_S(D)$ such that $H\cap (G_S(D)\setminus H)=\{v\}$ where $v$ is the vertex corresponding to $C$. It follows that $H$ consists of one or more blocks of $G_S(D)$. Similarly, the Seifert graph corresponding to the level 0 $s$-block also consists of one or more connected components. In the particular case that $G_S(D)$ is connected (in fact in this paper we only consider this case), the Seifert graph corresponding to the level 0 $s$-block is either a single vertex graph or a connected graph (which may contain more than one blocks).
}
\end{remark}

\medskip
From this point on, we will always assume that $D$ is a reduced alternating link diagram of some alternating link $L$. We shall use the notation $\G_S(D)$ for the simple graph obtained from $G_S(D)$ whose edges are assigned signed weights in the following way: the absolute value of the weight of an edge $e$ connecting two vertices $v_1$ and $v_2$ in $\G_S(D)$ is the total number of crossings between the $s$-circles $C_1$ and $C_2$ of $D$ corresponding to $v_1$ and $v_2$, and the sign of the weight of $e$ is the sign of these crossings. By Remark \ref{sign_remark}, all edges belonging to the same block in $\G_S(D)$ must have weights of the same sign. Since $D$ is reduced (that is, it is free of nugatory crossings), $G_S(D)$ is free of bridge edges hence $\G_S(D)$ is free of bridge edges whose absolute weights equal to one.

\begin{figure}[htb!]
\begin{center}
\includegraphics[scale=.45]{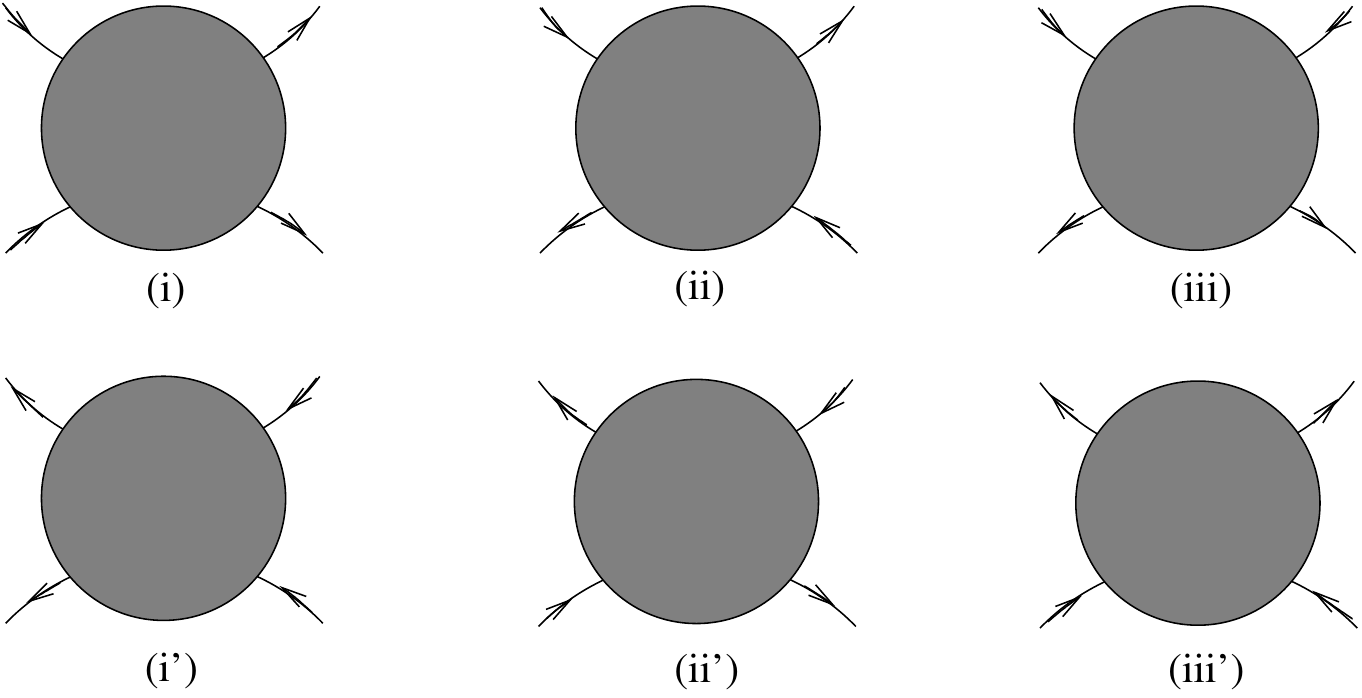}
\end{center}
\caption{\label{tangle_orientation} The possible combinations of the orientations of the strands entering/exiting $T$.
}
\end{figure}

\medskip
In the following we will examine how flypes on reduced alternating link diagrams can affect their Seifert graphs. Let us first consider the case when two reduced alternating link diagrams $D_1$ and $D_2$ are related by a single flype move. Let $T$ be the tangle where the flype takes place. Without loss of generality, let us assume that the crossing used for this flype move is on the right side of $T$ as shown in Figure \ref{flype}. There are six possible combinations to consider depending on the orientations of the four strands that enter/exit $T$. These are divided into three pairs (i) and (i'), (ii) and (ii'), and (iii) and (iii'), as shown in Figure \ref{tangle_orientation}. 

\medskip

\begin{figure}[htb!]
\begin{center}
\includegraphics[scale=.65]{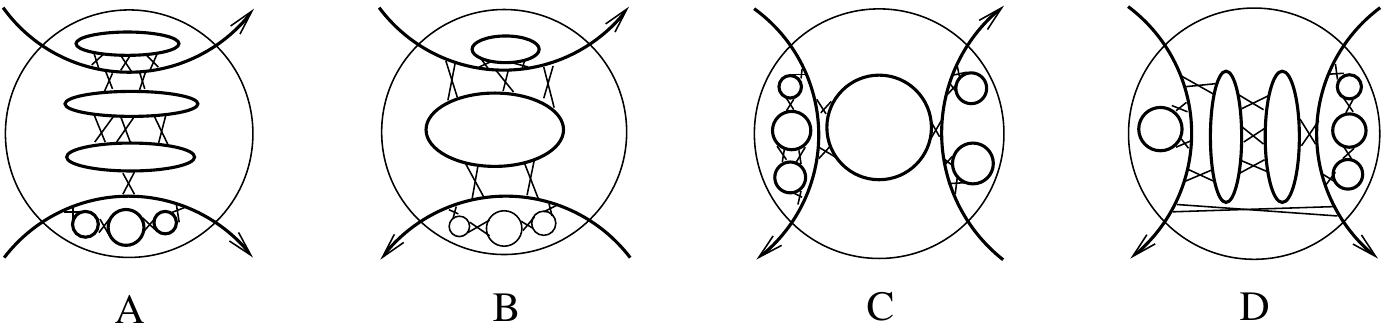}
\end{center}
\caption{\label{tangle_decomposition} The possible combinations of the arcs entering/exiting $T$ after all crossings are smoothed.
}
\end{figure}

By symmetry, we only need to consider one case from each pair, so we will only consider (i), (ii) and (iii). Furthermore, observe that in the Seifert circle decompositions $S(D_1)$, $S(D_2)$ of $D_1$, $D_2$, there are no crossings hence each strand entering $T$ must belong to a simple arc within $T$ which exits $T$ via an exiting strand, and the two arcs so formed do not cross each other. Depending on how these arcs are formed, the cases (i)--(iii) are further divided into the four sub-cases as illustrated in Figure \ref{tangle_decomposition}, which are denoted by A, B, C and D respectively. Some examples of Seifert circles and the smoothed crossings are included in Figure \ref{tangle_decomposition}. In the following we shall examine these sub-cases one by one and show that in each case $G_S(D_2)$ can be obtained from $G_S(D_1)$ with a suitable Whitney flip. Keep in mind that $S(D_1)$ and $S(D_2)$ are identical in the parts not involving the tangle $T$ and the crossing $O$.

\medskip

\medskip
In the cases below, the dashed lines in partial Seifert graphs of the bottom figures indicate the remaining parts of the graphs not associated with the structure $T \cup O$.

\begin{figure}[htb!]
\begin{center}
\includegraphics[scale=0.55]{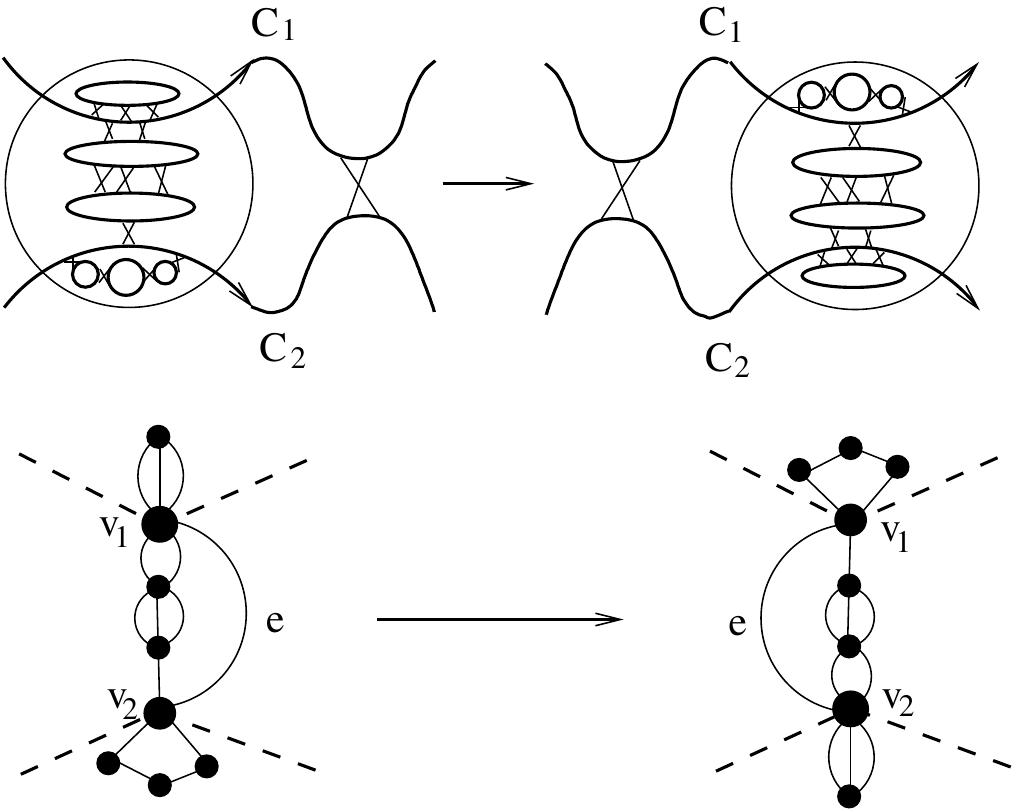}
\end{center}
\caption{Top: The effect of a flype move on $T$ and $O$ to the Seifert circle decompositions of $D_1$ and $D_2$ in the case that $T$ is of Type A; Bottom: The effect of the same operation on $G_S(D_1)$ and $S_S(D_2)$ is equivalent to a Type A Whitney flip as defined in Definition \ref{TypeAFlip}.
}
\label{TypeA}
\end{figure}

\medskip
\textbf{Case A.}  The top of Figure \ref{TypeA} illustrates how $S(D_1)$, 
$S(D_2)$ are related by the flype at $T$ and $O$. 

\medskip

Notice that in this case, the two Seifert circles $C_1$ and $C_2$ (defined by the two parallel arcs entering/exiting $T$ after all crossings are smoothed) must be distinct from each other (due to their orientations) and remain in both $S(D_1)$ and $S(D_2)$, while the crossing $O$ is moved to the left side of the tangle and the other structure in $T\cap S(D_2)$ is obtained from $T\cap S(D_1)$ by an 180$^\circ$ rotation around the horizontal centerline of $T$. At the Seifert graph level, the subgraph $G_1$ of $G_S(D_1)$ corresponding to the Seifert circles contained in $T$ is connected to the rest of graph through the vertices $v_1$, $v_2$ corresponding to $C_1$ and $C_2$. In other word, $v_1$ and $v_2$ form a 2-cut for $G_1$, and the subgraph $G_2$ of $G_S(D_2)$ corresponding to the Seifert circles contained in $T$ is obtained from $G_1$ by a Type A Whitney flip using $v_1$, $v_2$ as the 2-cut and the edge corresponding to the crossing $O$, as illustrated in the bottom of Figure \ref{TypeA}.

\begin{figure}[htb!]
\begin{center}
\includegraphics[scale=0.55]{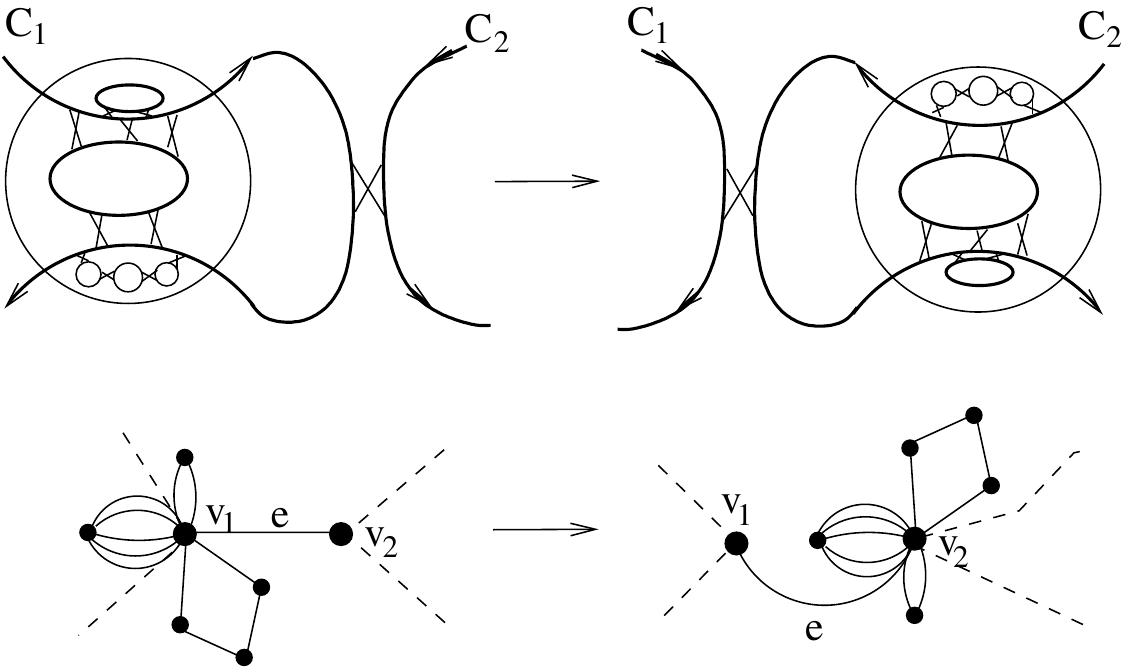}
\end{center}
\caption{Top: The effect of a flype move on $T$ and $O$ to the Seifert circle decompositions of $D_1$ and $D_2$ in the case that $T$ is of Type B; Bottom: The effect of the same operation on $G_S(D_1)$ and $S_S(D_2)$ is equivalent to a Type B Whitney flip as defined in Definition \ref{TypeBFlip}.
}
\label{TypeB}
\end{figure}

\medskip
\textbf{Case B.} The top of Figure \ref{TypeB} illustrates how $S(D_1)$, 
$S(D_2)$ are related by the flype at $T$ and $O$. 

In this case, the two Seifert circles $C_1$ and $C_2$ must also be distinct from each other and remain in both $S(D_1)$ and $S(D_2)$, while the crossing $O$ being a crossing between them both before and after the flype. At the Seifert graph level, the subgraph $G_1$ of $G_S(D_1)$ corresponding to the Seifert circles contained in $T$ is a union of blocks with $v_1$ being their defining cut vertex. $G_S(D_2)$ is obtained from $G_S(D_1)$ by performing a Type B Whitney flip on $G_1$ around $e$ as illustrated in the bottom of Figure \ref{TypeB}.

\begin{figure}[htb!]
\begin{center}
\includegraphics[scale=0.55]{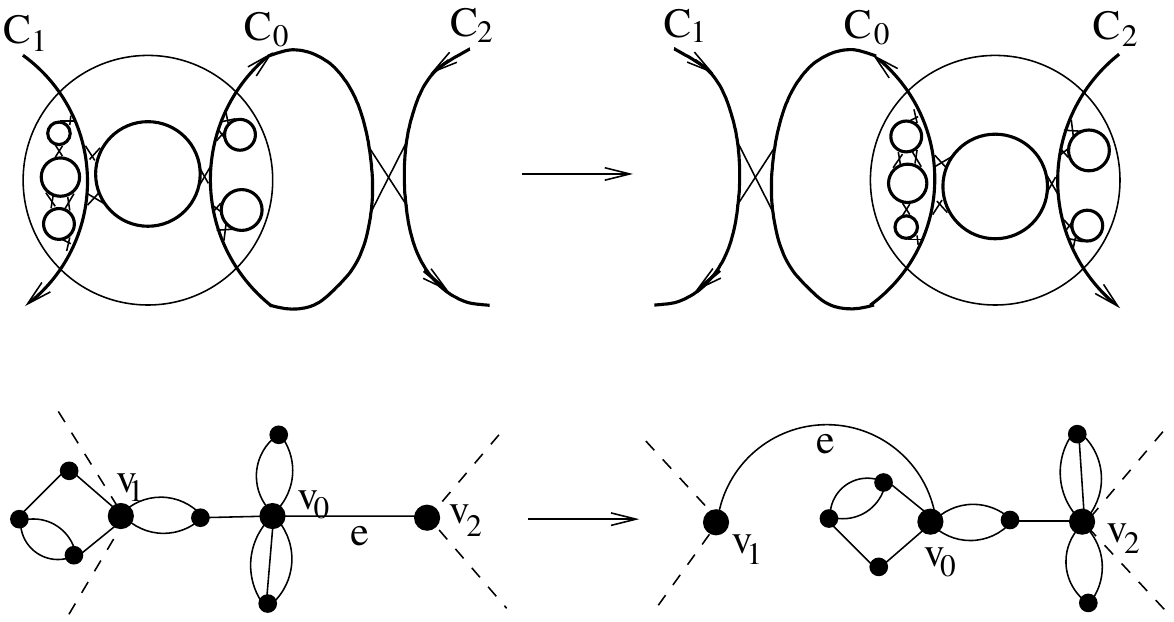}
\end{center}
\caption{Top: The effect of a flype move on $T$ and $O$ to the Seifert circle decompositions of $D_1$ and $D_2$ in the case that $T$ is of Type C; Bottom: The effect of the same operation on $G_S(D_1)$ and $S_S(D_2)$ is equivalent to a Type C Whitney flip as defined in Definition \ref{TypeCFlip}.
}
\label{TypeC}
\end{figure}

\medskip
\textbf{Case C.} The top of Figure \ref{TypeC} illustrates how $S(D_1)$, 
$S(D_2)$ are related by the flype at $T$ and $O$. 

In this case, the Seifert circles $C_1$ and $C_2$ must be distinct from each other and remain in both $S(D_1)$ and $S(D_2)$. The crossing $O$ is a single crossing between the Seifert circle $C_0$ and $C_2$ as shown in Figure  \ref{TypeC} in $D_1$, and becomes a single crossing between $C_0$ and $C_1$ in $D_2$.  At the Seifert graph level, the subgraph $G_1$ of $G_S(D_1)$ corresponding to the Seifert circles contained in $T$ is $v_1$-dependent in $G_S(D_1)\setminus e$ where $e$ is the unique edge between $v_0$ and $v_2$ corresponding to $O$. As shown by the bottom of Figure \ref{TypeC}, $G_S(D_2)$ is obtained from $G_S(D_1)$ by a Type C Whitney flip as defined in Definition \ref{TypeCFlip}.

\begin{figure}[htb!]
\begin{center}
\includegraphics[scale=0.5]{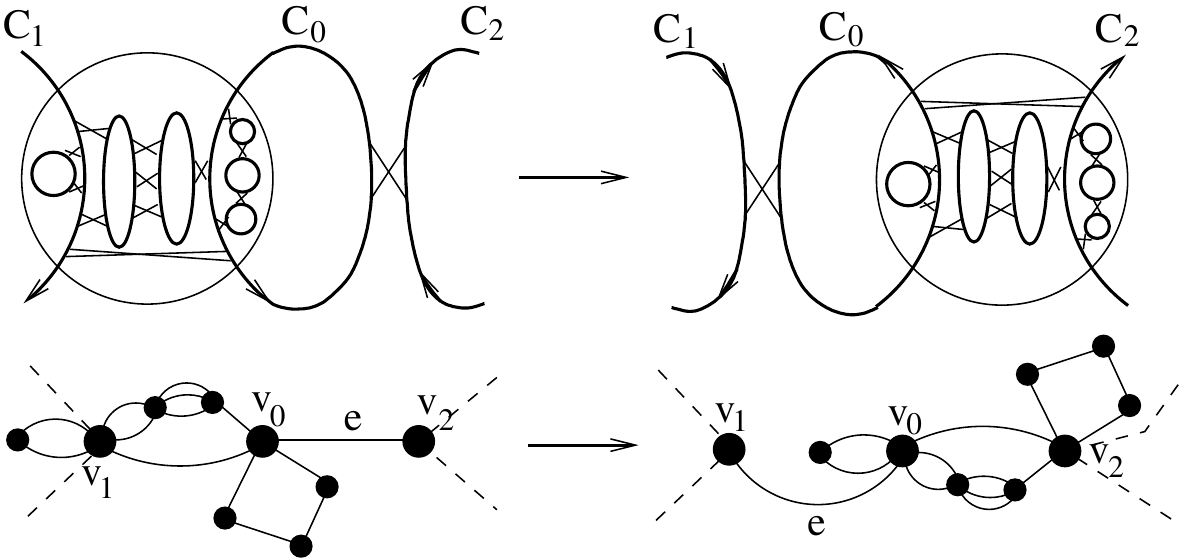}\quad \includegraphics[scale=0.5]{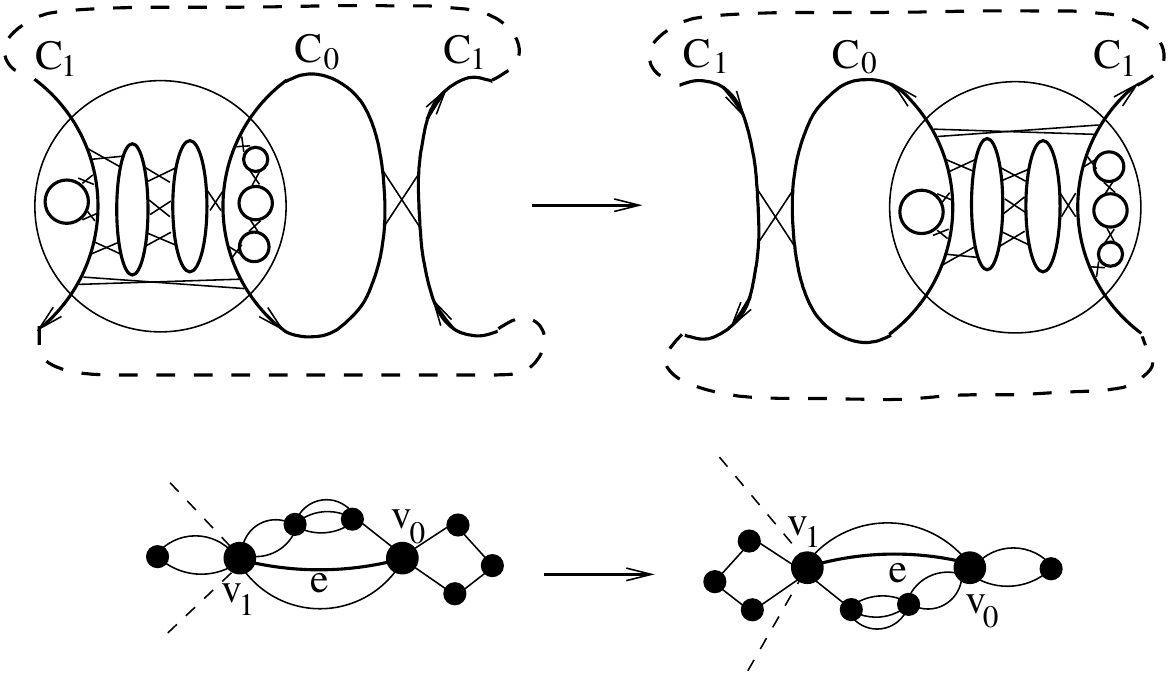}
\end{center}
\caption{Left: The effect of a flype move on $T$ and $O$ to the Seifert circle decompositions of $D_1$ and $D_2$ in the case that $T$ is of Type D and $C_1$ and $C_2$ are distinct and the effect of the same operation on $G_S(D_1)$ and $S_S(D_2)$; Right: The case when  $C_1=C_2$.
}
\label{TypeD}
\end{figure}

\medskip
\textbf{Case D.} There are two sub-cases for Case D. The left side of Figure \ref{TypeD} illustrates how $S(D_1)$,  $S(D_2)$ are related by the flype at $T$ and $O$ in the first sub-case, where the $s$-circles $C_1$ and $C_2$ are distinct from each other, and the crossing $O$ is a single crossing between the Seifert circle $C_0$ and $C_2$ as shown in Figure  \ref{TypeD} in $D_1$, and becomes a single crossing between $C_0$ and $C_1$ in $D_2$.  At the Seifert graph level, the subgraph $G_1$ of $G_S(D_1)$ corresponding to the Seifert circles contained in $T$ is $v_1$-dependent in $G_S(D_1)\setminus e$ where $e$ is the single edge between $v_0$ and $v_2$ corresponding to $O$. As shown in the bottom left of Figure \ref{TypeD} and similar to the case of Type C, $G_S(D_2)$ is obtained from $G_S(D_1)$ by a Type C Whitney flip as defined in Definition \ref{TypeCFlip}.

\begin{figure}[htb!]
\begin{center}
\includegraphics[scale=0.6]{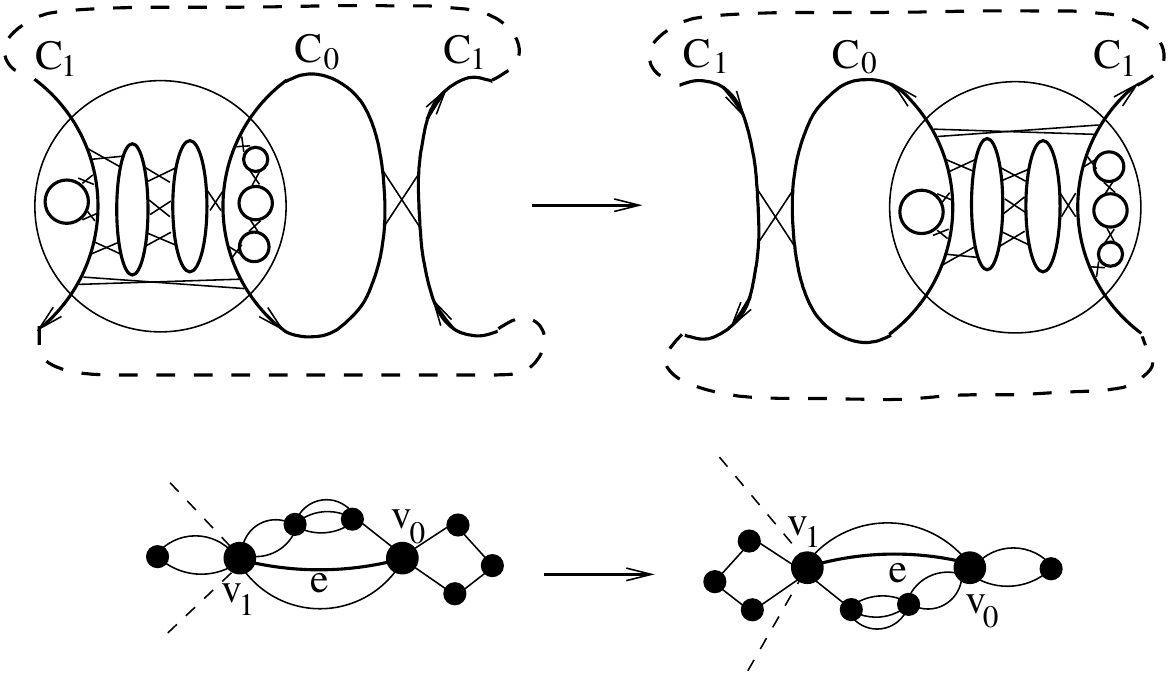}
\end{center}
\caption{Top: The effect of a flype move on $T$ and $O$ to the Seifert circle decompositions of $D_1$ and $D_2$ in the case that $T$ is of Type D and $C_1=C_2$ are the same Seifert circle; Bottom: The effect of the same operation on $G_S(D_1)$ and $S_S(D_2)$ is equivalent to a Type A Whitney flip as defined in Definition \ref{TypeAFlip}.
}
\label{TypeD2}
\end{figure}

\medskip
In the second sub-case we have $C_1=C_2$ as shown in the right side of Figure \ref{TypeD2}. In this case, the subgraph $G_2$ of $G_S(D_1)$ corresponding to the Seifert circles (and the crossings among them) contained in $T$ as well as the crossing $O$ is $v_1$-dependent and $G_S(D_2)$ is obtained from $G_S(D_1)$ by performing a special case of Type A Whitney flip as defined in Definition \ref{TypeAFlip}. Notice that this is special since $v_0$ is not connected to $v_1$ by any path not going through $G_2$, while this is allowed in the definition of a Type A Whitney flip.

\medskip

\medskip
The above analysis leads to the following theorem.

\medskip
\begin{theorem}\label{T1}
If two link diagrams $D_1$ and $D_2$ are related by a sequence of flype moves, then their Seifert graphs $G_S(D_1)$ and $G_S(D_2)$ are related by a sequence of Whitney flips of Types A, B or C, that is, $G_S(D_1)\stackrel{f}{\equiv} G_S(D_2)$.
\end{theorem}

\section{Writhe-like Invariants of Reduced Alternating Link Diagrams}\label{main}

Theorem \ref{T1} and Remark \ref {Whitney_flip_remark} allow us to derive various new writhe-like invariants from the Seifert graphs of reduced alternating link diagrams. The following is a list of such examples. 

\medskip
\begin{definition}\label{quantities}{\em
Let $D$ be a reduced alternating link diagram that is un-splittable ({\em i.e.} $G_S(D)$ is connected). Let $\{B_1, B_2, ..., B_k\}$ be the set of all blocks in $\G_S(D)$.

\medskip
\noindent
(i) $w_+(D)$ and $w_-(D)$: the sum of all positive edge weights and the sum of all negative edge weights in $\G_S(D)$ respectively;
\\

\smallskip
\noindent
(ii) $\xi_+(D)$ and $\xi_-(D)$: the number of positive blocks and the number of negative blocks in $\G_S(D)$ respectively, where a block is said to be {\em positive (negative)} if the weights in it are of positive (negative) sign;
\\

\noindent\smallskip
(iii) $W(D)$: the set of all edge weights in $\G_S(D)$;
\\

\noindent\smallskip
(iv) $w_B(D)=\{w(B_1), w(B_2), ..., w(B_k)\}$ where $w(B_j)$ is the sum of weights of edges of $B_j$;
\\

\noindent\smallskip
(v) $W_B(D)=\{W(B_1), W(B_2), ..., W(B_k)\}$ where $W(B_j)$ is  the set of all edge weights in $B_j$;
\\

\noindent\smallskip
(vi) $\beta(D)$: The Betti number $\beta(D)$ of $\G_S(D)$, which is defined as $\beta(D)=e(D)-v(D)+1$, where $e(D)$ is the number of edges in $\G_S(D)$ and $v(D)$ is the number of vertices in $\G_S(D)$;
\\

\noindent\smallskip
(vii) $\hat{\beta}(D)$: which is defined as the set $\{\beta(B_1), \beta(B_2), ..., \beta(B_k)\}$;
\\

\noindent\smallskip
(viii) $\Gamma(D)$: which is defined as the set $\{\gamma(B_1), \gamma(B_2), ..., \gamma(B_k)\}$ where $\gamma(B_j)$ is the length of a longest cycle in a block $B_j$ ($\gamma(B_j)=0$ if $B_j$ contains no cycles).
}
\end{definition}

\medskip
\begin{theorem}\label{T2}
Each quantity defined in Definition \ref{quantities} is a writhe-like invariant, that is, it is a link invariant within the space of all reduced alternating link diagrams.
\end{theorem}

\medskip
The proof of Theorem \ref{T2} is straight forward. For example,  $\beta(D)$ and $\Gamma(D)$ are writhe-like invariants since a Whitney flip cannot change a cycle except possibly the position of an edge with absolute weight one in the cycle. We leave the details to the reader. 

\medskip
We say that an invariant is {\em stronger} than another one if it can distinguish all link types the other invariant can, but not vice versa. For example, $w_+(D)$ and $w_-(D)$ together is stronger than $w(D)$, $W(D)$ is stronger than $w_+(D)$ and $w_-(D)$ combined, $w_B(D)$ is stronger than $\xi_+(D)$ and $\xi_-(D)$ combined, and $W_B(D)$ is stronger than $w_B(D)$ and also stronger than $W(D)$. Of course, some of the invariants in the list are not stronger to each other. For example, there are distinct reduced alternating link diagrams that $\beta(D)$ can distinguish, but not $w_B(D)$, and vice versa. 

\medskip
With some effort, one can expand the list of writhe-like invariants in Definition \ref{quantities} that use some  diagrammatic features of $D$ itself. We demonstrate below one such example.

\begin{figure}[htb!]
\begin{center}
\includegraphics[scale=1.1]{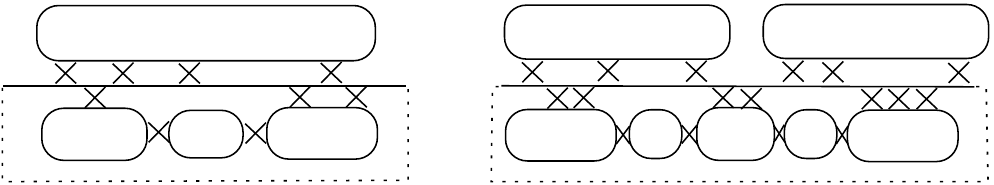}
\end{center}
\caption{Left: A diagram containing a locked pair of Seifert circles with 4 clusters; Right: A diagram containing no locked Seifert circle pairs.
\label{lock}
}
\end{figure}

\medskip
Consider two vertices $v_1$, $v_2$ of $\G_S(D)$ that are attached to the same cut vertex $v$ by edges of different sign (so $v$ is necessarily a cut vertex). Let $C$, $C_1$ and $C_2$ be the Seifert circles corresponding to $v$, $v_1$ and $v_2$ respectively. Let $c(C_1)$ and $c(C_2)$ be the sets of crossings between $C$ and $C_1$, $C$ and $C_2$ respectively. If we travel along $C$ once, we will encounter each crossing in $c(C_1)$ and $c(C_2)$ exactly once. A set of crossings from $c(C_1)$ ($c(C_2)$) is said to form a {\em cluster} if it is a maximal set with the property that we can use one of them as a starting point to travel along $C$ to reach the last crossing in the set without encountering any crossing in $c(C_2)$ ($c(C_1)$). We say that $v_1$ and $v_2$ are {\em locked} if the total number of clusters is 4 or more. The left side of Figure \ref{lock} shows (a portion of) a link diagram containing a locked pair with Seifert circles with 4 clusters and the right side is an example containing no locked Seifert circle pairs.

\medskip

\medskip
Let us consider what happens to a locked pair of Seifert circles when a flype is applied. We wish to show that a flype cannot break the lock nor change the number of clusters. If the flype is of Type A or Type B,  then all crossings between the two Seifert circles are either all within or outside the tangle except that one of them may be used as  the crossing for the flype, but either way the flype does not change the clusters. If the flype is of Type C, then the crossing used for the flype cannot be a crossing in the lock since then the number of clusters can only be 2 between $C_1$ and $C_2$ if one of them has a single crossing with $C_0$. Thus all crossings between $C_1$ and $C_2$ are either within the tangle or outside of the tangle, either way the flype will not change the clusters. If the flype is of Type D, then the first case as shown in the left side of Figure \ref{TypeD} is similar to the case of a Type C flype. If it is the second case as shown in the fight side of Figure \ref{TypeD} (where $C_1=C_2$), so the Seifert circle locked with $C_1$ must be within $C_0$, and the effect of the flype to the pair is merely a reflection hence the number of clusters does not change either. Thus a flype will not change the clusters of a lock, which means the number of crossings in each cluster and the cyclic order of these clusters along the Seifert circle $C_0$ used in the definition of the lock are all preserved. These observations lead to the following theorem.

\medskip
\begin{theorem}\label{Lock_Theo}
Define $\Phi(D)=\{\vec{\phi}_1, \vec{\phi}_2,...,\vec{\phi}_m\}$ where $m$ is the total number of locked pairs of Seifert circles and  $\vec{\phi}_j$  is a vector consisting of the number of clusters in the $j$-th lock according to their cyclic order along the orientation of the Seifert circle used to define the lock.
\end{theorem}

\section{Applications and Examples}\label{applications}

While the writhe-like invariants listed in Definition \ref{quantities} may not seem to be strong link invariants, the combined use of them can in fact be quite powerful. For example, we can use them to distinguish most alternating knots with small crossing numbers, including all alternating knots up to 9 crossings. Example \ref{E1} shows that some of these invariants can distinguish two large alternating knots when other invariants such as the Jones polynomial is difficult to compute. 

\medskip
\begin{example}\label{E1}{\em
Figure \ref{mutation_graph} shows the Seifert circle decompositions of two reduced alternating knot diagrams $D_1$ and $D_2$ with the crossings in the middle being negative. Both have 19 crossings. $D_1$ and $D_2$ are in the closed braid form so $\G_S(D_1)$ and $\G_S(D_1)$ are identical: both are a simple path consisting of three consecutive edges of weights $6$, $-5$ and $8$ respectively. However we have $\Phi(D_1)=\{(2,-1,3,-1,1,-3),(1,-1,3,-1,4,-3)\}$ and $\Phi(D_2)=\{(2,-1,3,-1,1,-3),(3,-1,3,-1,2,-3)\}$. Thus $D_1$ and $D_2$ are different knots since the second vectors in $\Phi(D_1)$ and $\Phi(D_2)$ cannot be obtained from each other by a cyclic order permutation.

\medskip
\begin{figure}[htb!]
 \begin{center}
\includegraphics[scale=0.4]{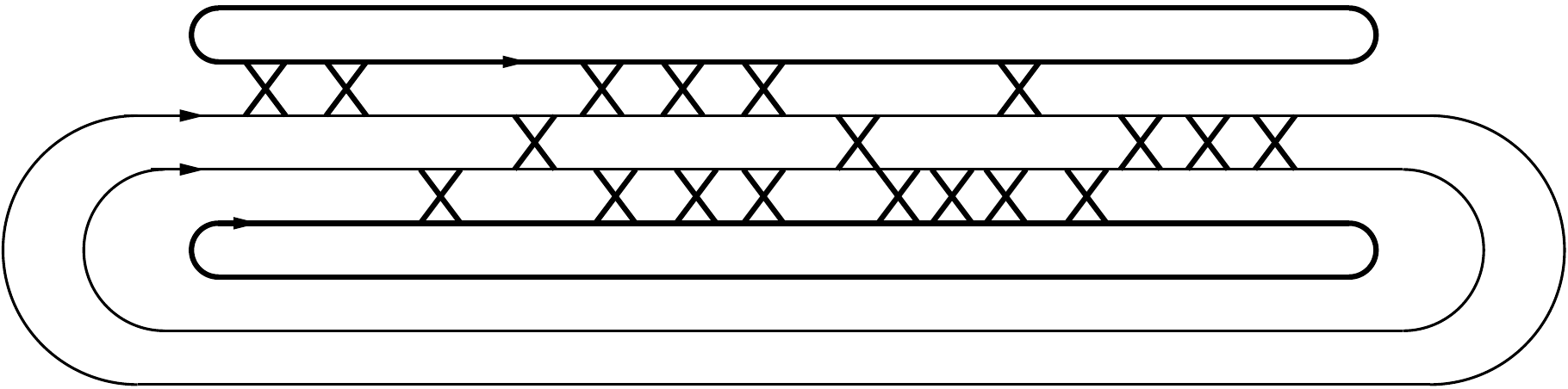}\\
\includegraphics[scale=0.4]{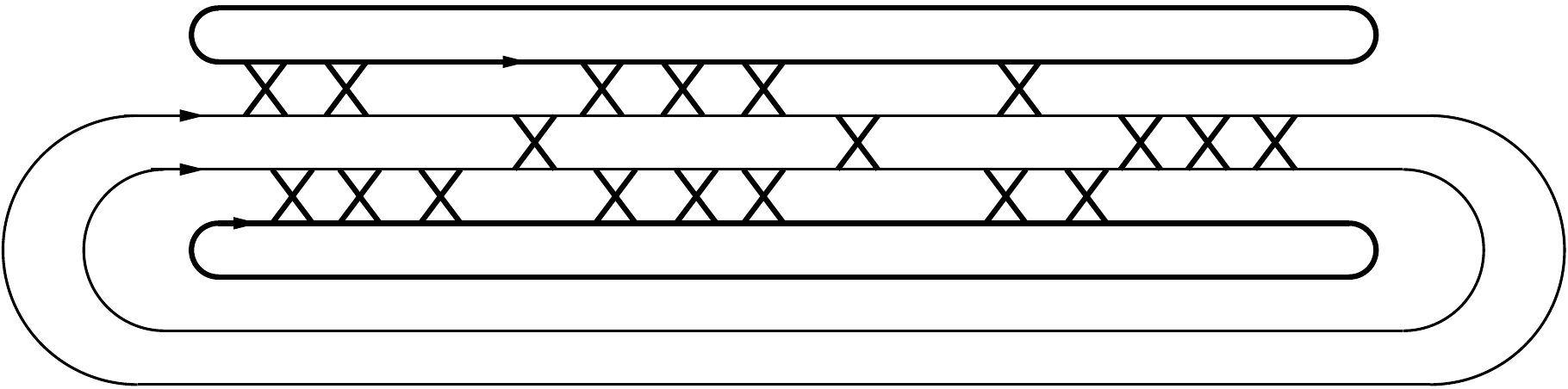}
\end{center}
\caption{The Seifert circle decompositions of two reduced alternating knot diagrams $D_1$ (top) and $D_2$ (bottom). The signs of the crossings are indicated in the figures.
\label{mutation_graph}
}
\end{figure}
}
\end{example}

\begin{figure}[htb!]
 \begin{center}
\includegraphics[scale=0.5]{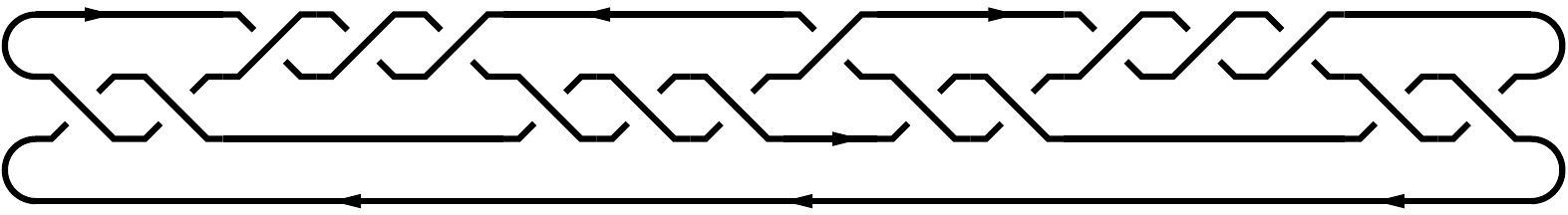}
\end{center}
\caption{The reduced alternating diagram $D$ of the rational knot corresponding to the rational number $278/641$ which has the continued fraction decomposition vector $[2,3,3,1,2,3,2]$.
\label{mutation}
}
\end{figure}

\medskip
Our next example is the proof of a theorem that contains a new result concerning oriented rational links. Thus we need some preparations before we can state and prove the theorem. We shall assume that our readers have some basic knowledge about rational links. \cite{Burde} is a good reference for readers who are not familiar with these subjects. Let $p$ and $q$ be two positive integers with $\rm{gcd}(p,q)=1$ and $0<p<q$. Let $(a_1,a_2,...,a_{2k+1})$ be the (unique) vector of odd length with positive integer entries such that 
$$
\displaystyle\frac{p}{q}=\displaystyle\frac{1}{a_1+\frac{1}{a_2+\frac{1}{....\frac{1}{a_{2k}+\frac{1}{a_{2k+1}}}}}}.
$$
For the sake of convenience we will write the above as $p/q=[a_1,a_2,...,a_{2k+1}]$. It is known that any rational link has a reduced alternating diagram called a 4-plat corresponding to the odd length continued fraction decomposition vector of some positive integers with $\rm{gcd}(p,q)=1$ and $0<p<q$ as shown in the Figure \ref{mutation} for the case of $p/q=278/641=[2,3,3,1,2,3,2]$. We will write this rational link as $L(p/q)$.

\medskip

\medskip
It is known that rational links are invertible, that is, changing the orientations of all components in a rational link will not change the link type. However, in the case that an oriented rational link $L$ has two components (this happens if and only if $q$ is even if $L$ is represented by a 4-plat corresponding to $p/q=[a_1,a_2,...,a_{2k+1}]$), changing the orientation of only one component of $L$ may result in a rational link that is topologically different from $L$ (as oriented links). $L$ is said to be {\em strongly invertible} if changing the orientation of one of its component does not change its link type. \cite[Theorem 8.1]{Ka-La} states that if a two-component rational link $L(p/q)$ is strongly invertible, then $p/q=[a_1,a_2,...,a_k,\alpha,a_k,...,a_2,a_1]$ for some integers $a_1>0$, ..., $a_k>0$, $\alpha> 0$. The  following theorem  strengthens this result.

\medskip
\begin{theorem}\label{T3}
Let $L(p/q)$ be a two-component oriented rational link (so that $\rm{gcd}(p,q)=1$, $0<p<q$ and $q$ is even). Then 
 $L(p/q)$ is strongly invertible if and only if $p/q=[a_1, a_2, ..., a_n, \alpha, a_n, ..., a_2, a_1]$ where $a_1, a_2, ..., a_n$ are positive integers and $\alpha>0$ is odd.
\end{theorem}

\medskip
\begin{proof} Without loss of generality, the bottom long arc of $L(p/q)$ in its 4-plat form as shown in the top of Figure \ref{mutation}. Let $L_1(p/q)$ and $L_2(p/q)$ denote the two rational links by assigning the second component different orientation so that the first crossing from the left is positive in $L_1(p/q)$ and negative in $L_2(p/q)$.

\medskip
Let us first examine the case when $L(p/q)$ (with $q$ even) is strongly invertible, that is, $L_1(p/q)\sim L_2(p/q)$ and $L(p/q)$ is either $L_1(p/q)$ or $L_2(p/q)$. By \cite[Theorem 8.1]{Ka-La}, $p/q=[a_1,a_2,...,a_k,\alpha,a_k,...,a_2,a_1]$ for some integers $a_1>0$, ..., $a_k>0$, $\alpha> 0$. We claim that in this case the right most crossing is negative (positive) in $L_1(p/q)$ ($L_2(p/q)$). If this is not true, say that the right most crossing in $L_1(p/q)$ is also positive, then as we travel along the $s$-circle containing the long arc starting from the left most crossing, the first and last blocks of $L_1(p/q)$ we encounter are both positive, which implies that $\xi_+(L_1(p/q))=\xi_-(L_1(p/q))+1$. On the other hand, for $L_2(p/q)$ the first block we encounter is negative, which implies that $\xi_-(L_2(p/q))\ge \xi_+(L_1(p/q))$, which is a contradiction. This implies that the strands at the left and right ends of the 4-plat $L_1(p/q)$ are as shown in Figure \ref{fig_strong1}, in which the middle points marked by 1, 2, 3 and $1^\p$, $2^\p$, $3^\p$ are the points to be connected to create the middle section of $L_1(p/q)$ corresponding to the $\alpha$ crossings. If we start at the point marked by 3 from the left and travel according to the orientation of the strand, we will end at one of the points marked by 1, 2 or 3 in the middle. Say we stop at 3 (the other cases can be similarly discussed and are left to the reader). Then  since the continued fraction decomposition vector of $p/q$ is a palindrome, if we start at the point marked by $3^\p$ from the right and travel according to the orientation of the strand, we will also end at the point marked by $3^\p$ in the middle. If we continue, we will now have to go through the crossings in the middle, since otherwise we will have to connect 3 to $3^\p$ by a straight line segment but that will violate the orientations. For the same reason, $\alpha$ has to be odd since otherwise we will still have to connect 3 to $3^\p$ by a strand. Hence the only possibility is that $\alpha$ is odd and the point 3 is connected to the point $2^\p$ by going through these crossings.

\medskip
\begin{figure}[!hbt]
\begin{center}
\includegraphics[scale=0.6]{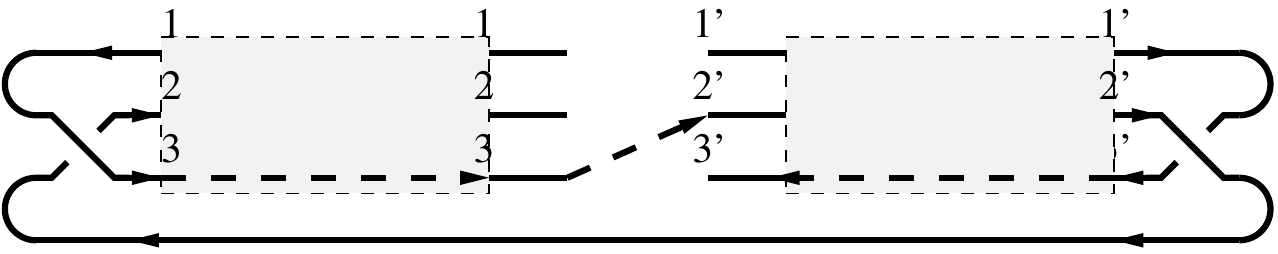}
\end{center}
\caption{The orientations of the strands at the ends of a strongly invertible rational link in a 4-plat.
\label{fig_strong1}}
\end{figure}

\medskip
Now, let us consider the case when $p/q=[a_1, a_2, ..., a_n, \alpha, a_n, ..., a_2, a_1]$ where $a_1, a_2, ..., a_n$  are positive integers and $\alpha>0$ is odd, and $q$ is even (so $L(p/q)$ has two components). Similar to the discussion in the above, let us consider the partial 4-plat of $L_1(p/q)$ as shown in Figure \ref{fig_strong2}. Notice in this case we do not know the orientation of the strand at the right end of the 4-plat that does not belong to the bottom long strand. 

\medskip
\begin{figure}[!hbt]
\begin{center}
\includegraphics[scale=.6]{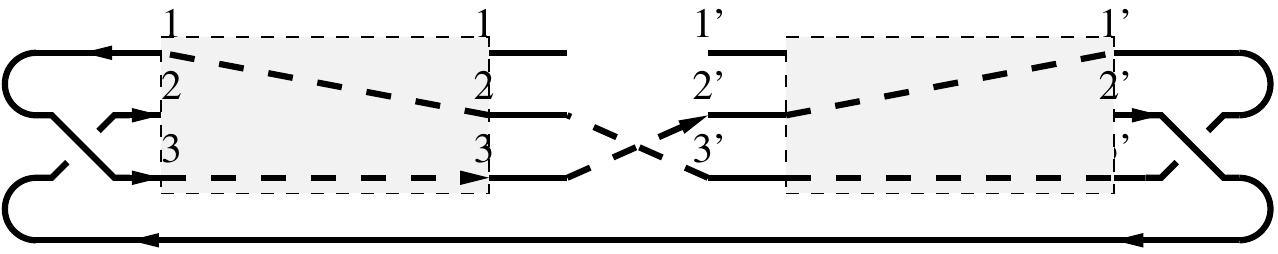}\\

\bigskip
\includegraphics[scale=.6]{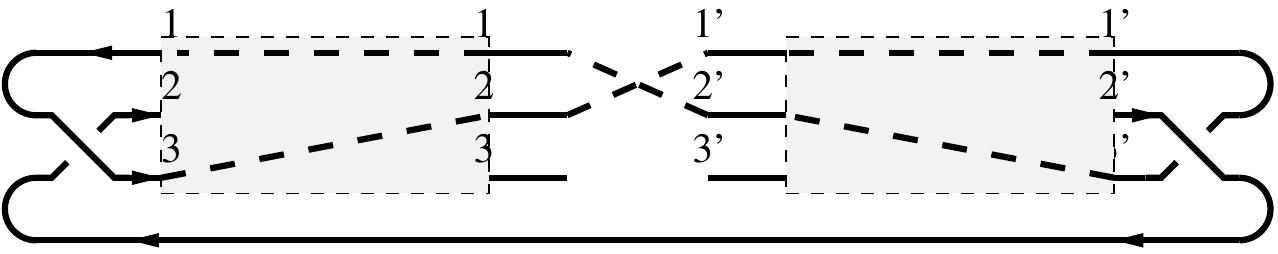}
\end{center}
\caption{How to strands connect in a rational link in a symmetric 4-plat form.
\label{fig_strong2}}
\end{figure}

\medskip
If we start at the point marked by 3 from the left and travel according to the orientation of the strand, then again we will end at one of the points marked by 1, 2 or 3 in the middle. Say we stop at 3. 3 cannot be connected to $3^\p$ by going through the crossings in the middle since $\alpha$ is odd. So if 3 is connected to $3^\p$ then it is by a straight line segment and the middle crossings will be between points 1, 2 and $1^\p$, $2^\p$ in the middle. Since the link has two components, it is necessary that 2 be connected to $2^\p$ and 1 be connected to $1^\p$. This is impossible since $\alpha $ is odd. Thus 3 can only be connected to $2^\p$ by going through the $\alpha$ crossings in the middle, and we shall travel to $1^\p$ to the right side in order for the link to have two components. This then determines the orientation of the strand at the right side (that does not belong to the bottom long strand), which is the same as that of Figure \ref{fig_strong1}. If we now take a $180^\circ$ rotation of the 4-plat around the $y$-axis, then we obtain $L_2(p/q)$ (it uses the same vector $[a_1, a_2, ..., a_n, \alpha, a_n, ..., a_2, a_1]$ and starts with a negative crossing on the left). 

\medskip
There are two cases left. If we start at 3 and end at 1 in the middle, then the discussion is identical to the above. If we start at 3 and end at 2 in the middle, then we have to go through the middle crossings to reach either $1^\p$ or $3^\p$ as shown in the bottom of Figure \ref{fig_strong2} for one of the two cases. The remaining argument is similar to the above. \qed
\end{proof}

\section*{Acknowledgement} The authors thank Dr. Christoph Lamm for pointing out an error in an earlier version of this manualsctipt.


\bigskip

\section*{References}

\begin{thebibliography}{99}
\bibitem{Adams} C.~Adams,
{\em The Knot Book}, American Mathematical Soc., 2004.

\bibitem{Burde} G.~Burde and H.~Zieschang, {\em Knots}, Walter de Gruyter, 2003.

\bibitem{DHL2019} Y. Diao, G. Hetyei and P. Liu. The Braid Index Of Reduced Alternating Links,  {\em Math. Proc. Cambridge Phil. Soc.}, doi:10.1017/S0305004118000907, 1--19 (January 2019).

\bibitem{Ka-La}
L.H.~ Kauffman, and S.~ Lambropoulou, 
{\em On the classification of rational knots,}
Enseign. Math \textbf{2} (2002), 357--410.   

\bibitem{MT1} W.~Menasco and M.~Thistlethwaite, The Tait Flyping Conjecture. {\em Bull. Amer. Math. Soc.} \textbf{25} (1991), 403--412.

\bibitem{MT2} W.~Menasco and M.~Thistlethwaite, The Classification of Alternating Links. {\em Ann. Math.} \textbf{138} (1993), 113--171. 

\bibitem{Rolfsen} D.~Rolfsen, {\em Knots and Links}, AMS Chelsea Pub., 2003.

\bibitem{West} D.~West, {\em Introduction to Graph Theory} (Second edition), Prentice Hall 2001.

\bibitem{Whitney} H.~Whitney,  2-Isomorphic Graphs. {\em Amer. J. Math.} \textbf{55} (1-4) (1933), 245--254.

\end{thebibliography}
\end{document}